\documentclass[11pt,leqno]{article}

\usepackage{amsmath,amsfonts,amscd,amssymb,theorem}

\long\def\comment#1\endcomment{}

\comment
\endcomment

\makeatletter
\begingroup
\gdef\th@dotted{\normalfont\itshape
  \def\@begintheorem##1##2{%
        \item[\hskip\labelsep \theorem@headerfont ##1\ ##2.]}%
\def\@opargbegintheorem##1##2##3{%
   \item[\hskip\labelsep \theorem@headerfont ##1\ ##2\ (##3).]}}
\endgroup
\makeatother

\theoremstyle{dotted}

\newtheorem{theorem}{Theorem}[section]
\newtheorem{lemma}[theorem]{Lemma}

\makeatletter
\begingroup
\gdef\th@upshape{\normalfont
  \def\@begintheorem##1##2{%
        \item[\hskip\labelsep \theorem@headerfont ##1\ ##2.]}%
\def\@opargbegintheorem##1##2##3{%
   \item[\hskip\labelsep \theorem@headerfont ##1\ ##2\ (##3).]}}
\endgroup
\makeatother

\theoremstyle{upshape}

\newtheorem{defn}[theorem]{Definition}
\newtheorem{remark}[theorem]{Remark}

\makeatletter
\renewcommand{\subsection}{\@startsection{subsection}{2}{0pt}{-3ex
plus -1ex minus -0.2ex}{-2mm plus -0pt minus
-2pt}{\normalfont\bfseries}} \makeatother

\makeatletter
\@addtoreset{equation}{section}
\makeatother

\newcommand{\cntrct}                
{\hspace{2pt}\raisebox{1pt}{\text{$\lrcorner$}}\hspace{2pt}}

\newcommand{\proof}[1][Proof.]{\smallskip\noindent{\em #1}}
\def\endproof{\hfill\ensuremath{\square}\par\medskip}

\def\eqref#1{\thetag{\ref{#1}}}

\let\latexref=\ref
\def\ref#1{{\normalfont{\latexref{#1}}}}

\newcommand{\wt}{\widetilde}
\newcommand{\wh}{\widehat}

\newcommand{\idot}{{\:\raisebox{1pt}{\text{\circle*{1.5}}}}}
\newcommand{\hdot}{{\:\raisebox{3pt}{\text{\circle*{1.5}}}}}
\newcommand{\C}{{\mathbb C}}

\newcommand{\B}{{\cal B}}
\newcommand{\A}{{\cal A}}
\newcommand{\D}{{\cal D}}

\newcommand{\CC}{{\cal C}}

\newcommand{\cchar}{\operatorname{\sf char}}

\newcommand{\I}{{\cal I}}

\newcommand{\calo}{{\cal O}}

\newcommand{\h}{{\mathfrak h}}

\newcommand{\eps}{\varepsilon}
\renewcommand{\phi}{\varphi}

\newcommand{\bimod}{\operatorname{\sf bimod}}
\newcommand{\gr}{\operatorname{\sf gr}}

\newcommand{\Coh}{\operatorname{Coh}}
\newcommand{\Hom}{\operatorname{Hom}}
\newcommand{\Ext}{\operatorname{Ext}}

\newcommand{\RHom}{\operatorname{RHom}}
\newcommand{\Rhom}{\operatorname{{\bf R}{\cal H}{\it om}}}

\newcommand{\HH}{\operatorname{{\mathbb H}}}
\newcommand{\hh}{\operatorname{{\cal H}{\cal H}}}

\newcommand{\Ker}{\operatorname{{\sf Ker}}}
\renewcommand{\Im}{\operatorname{{\sf Im}}}

\title{Some remarks on formality in families}

\author{D. Kaledin\thanks{Partially supported by CRDF grant
RUM1-2694-MO05.}}

\begin{document}

\date{{\em To my first math teacher Victor Ginzburg, on the occasion
of his 50-th birthday}}

\maketitle

\tableofcontents

\section*{Introduction}

The notion of a {\em formal\/} differential graded algebra -- that
is, a DG algebra $A^\hdot$ quasiisomorphic to its cohomology algebra
$H^\hdot(A^\hdot)$ -- is by now a familiar sight in many areas of
mathematics; we can quote, for instance, the classic paper
\cite{SDGM}, where formality was established for the de Rham
cohomology algebra of a compact K\"ahler manifold $X$, which had
numerous applications to the topology of compact K\"ahler
manifolds. A well-known series of obstructions to formality is given
by the so-called {\em Massey products}. It would be very convenient
to know that these give the only obstruction -- a DG algebra with
Massey products is formal. Unfortunately, this is not true (for a
counterexample, see e.g. \cite{S}). Therefore in works such as
\cite{SDGM} Massey products play only a marginal role, and the main
technical tool is the notion of a {\em minimal model} introduced by
D. Sullivan.  However, this brings about some problems, of which the
most obvious one is that minimal models usually do not exist for
families of DG algebras over a sufficiently non-trivial base.

In this paper, we construct a certain refinement of the Massey
products which does characterise formality uniquely, and moreover,
behaves well for families of DG alebras. As an application, we prove
two results on formality in families.

Of course, at least morally, and at least in some cases, both
results are not new. However, it seems that accurate and complete
proofs are not available in the existing literature, which precludes
applications in non-standard setting. The goal of this paper is to
provide such a proof. To save space, we only sketch those proofs
that deal with DG algebras over a field, -- this material is quite
standard, -- and coversely, we try to be really precise when it
comes to families of DG algebras over a base. Our approach to
formality is motivated by and partially follows the paper \cite{H}:
we treat formality of a DG algebra $A^\hdot$ as triviality of the
normal cone deformation associated to the canonical filtration on
$A^\hdot$, and we use deformation theory methods to find criteria
for this triviality.

\subsection*{Acknowledgements.} This paper owes its existence to
M. Lehn: he convinced me that this is the case where writing down
detailed proofs is a meaningful thing to do, and read through
innumerable first drafts, each one longer than the preceding one
(and hopefully, less incorrect). I am also grateful to M. Verbitsky
for reading some of those drafts, giving me a second opinion, and
suggesting the reference \cite{S}.

\section{Kodaira-Spencer classes.}\label{fo}

Let $A$ be an associative, not necessarily commutative algebra. The
embedding of the diagonal defines a canonical map $A \otimes A \to
A$ of $A$-bimodules. Denote its kernel by $\I_A \subset A \otimes
A$. For any $A$-bimodule $M$, the groups
\begin{equation}\label{HH.defn}
HH^{\hdot+1}(A,M) = \Ext^\hdot_{A-\bimod}(\I_A,M)
\end{equation}
are called the (reduced) Hochschild cohomology groups of the algebra
$A$ with coefficients in $M$. (One can show that alternatively,
$HH^\hdot(A,M)=\Ext^\hdot_{A-\bimod}(A,M)$, which explaines the
shift in the index; we will not need this.) If $M=A$ is $A$ itself
considered as an $A$-bimodule, than the groups $HH^\hdot(A,A)$ are
denoted simply by $HH^\hdot(A)$.

A square-zero extension $\wt{A}$ of the algebra $A$ by the bimodule
$M$ is by definition an associative algebra $\wt{A}$ equipped with a
two-sided ideal $N \subset \wt{A}$ such that we are given an
isomorphism $\wt{A}/N \cong A$, the $\wt{A}$-bimodule structure on
$N$ factors through an $A$-bimodule structure, and we are given an
$A$-bimodule isomorphism $N \cong M$. Every square-zero extension
$\wt{A}$ defines a Hochschild cohomology class $\theta_{\wt{A}} \in
HH^2(A,M)$ by means of the following procedure. Denote by $\I_M
\subset M \otimes A$ the kernel of the natural $A$-bimodule map $M
\otimes A \to M$. Consider the $\wt{A}$-$A$-bimodule $\wt{A} \otimes
A$ ($\wt{A}$ acts by left multiplication, $A$ acts by right
multiplication). We have a short exact sequence of
$\wt{A}$-$A$-bimodules
$$
\begin{CD}
0 @>>> M \otimes A @>>> \wt{A} \otimes A @>>> A
\otimes A @>>> 0.
\end{CD}
$$
In particular, we have an embedding $\I_M \to \wt{A} \otimes A$
and a surjection $\wt{A} \otimes A \to A$. The middle cohomology
$\wt{\I}_{A,M}$ of the complex $\I_M \to \wt{A} \otimes A \to A$ appears
as the middle term of a short exact sequence
\begin{equation}\label{exx}
\begin{CD}
0 @>>> M @>>> \wt{\I}_{A,M} @>>> \I_A @>>> 0
\end{CD}
\end{equation}
of $\wt{A}$-$A$-bimodules. One checks easily that this sequence is
in fact a sequence of $A$-bimodules. We take $\theta_{\wt{A}} \in
HH^2(A,M)$ to be the Yoneda class of the extension \eqref{exx}.

A first-order deformation $\wt{A}$ of the algebra $A$ is by
definition a square-zero extension of $A$ by $A$; equivalently, it is
an associative algebra $\wt{A}$ equipped with a
$\wt{A}$-bimodule map $\eps:\wt{A} \to \wt{A}$ such that $\Ker\eps =
\Im\eps \subset \wt{A}$, and an algebra isomorphism $\wt{A}/\Ker\eps
\cong A$. Any first-order deformation $\wt{A}$ defines a class
$\theta_{\wt{A}} \in HH^2(A) = HH^2(A,A)$.

\medskip

There are many ways to present this construction. The one we have
chosen has the following advantage: it works without any changes for
a flat algebra $A$ in an arbitrary abelian tensor category $\CC$.

As a first application of this additional degree of freedom, we show
that the same definition can be used to study higher-order
deformations. Namely, by a formal deformation $\wt{A}$ of an algebra
$A$ in a symmetric tensor category $\CC$ we will understand an
associative algebra $\wt{A}$ equipped with an injective algebra map
$h:\wt{A} \to \wt{A}$ and an isomorphism $\wt{A}/h(\wt{A}) \cong
A$. A formal deformation $\wt{A}$ is an algebra in the tensor
category $\CC[h]$ of objects in $\CC$ equipped with an endomorphism
$h$; since $h:\wt{A} \to \wt{A}$ is injective, $\wt{A}$ is flat in
$\CC[h]$ if $A$ is flat in $\CC$. Given such a formal deformation
$\wt{A}$, we can consider a trivial first-order deformation
$\overline{A} = \wt{A}\langle\eps\rangle = \wt{A} \oplus
\wt{A}\cdot\eps$ of the algebra $\wt{A}$ in $\CC[h]$, and redefine
the endomorphism $h:\overline{A} \to \overline{A}$ by setting
$h_{new} = h_{old} + \eps$. Since $\eps^2=0$ and $h_{old}$ is
injective, $h_{new}$ is also injective -- $\Ker h_{new} \cap \Ker
\eps \subset \Ker \h_{old}$ must be trivial, by induction, $\Ker
h_{new} \cap \Ker \eps^l$ is then trivial for every $l \geq 1$, but
already $\Ker h_{new} \cap \Ker \eps^2 = \Ker h_{new}$. Thus algebra
$\overline{A}$ with the new endomorphism $h$ is still a first-order
deformation of the algebra $\wt{A}$ in $\CC[h]$. If $\wt{A} \cong
A[h]$ were a trivial formal deformation of $A$, then this first
order deformation is trivial; in general, however, it might be
non-trivial and defines a cohomology class
$$
\Theta_{\wt{A}} \in HH^2(\wt{A})
$$
called the {\em Kodaira-Spencer class} of the deformation
$\wt{A}$. This describes the non-triviality of the deformation
$\wt{A}$ (not completely -- we do not claim that in the general
case, $\Theta_{\wt{A}}=0$ implies $\wt{A} \cong A[[h]]$).

We note that the groups $HH^\hdot(\wt{A})$ are equipped with a
natural endomorphism $h$; we have a natural map $HH^2(\wt{A})/h \to
HH^2(A)$, and the image of the class $\Theta_{\wt{A}}$ under this
map is the cohomology class corresponding to the first-order
deformation $\wt{A}/h^2(\wt{A})$ of the algebra $A$.

\section{Explicit cocycles.}

Let us compare the formalism of Section~\ref{fo} with the more
standard approaches to the Kodaira-Spencer class. Firstly, assume
that the $A$-bimodule $M$ is injective as an object of the category
$\CC$ (for example, this is always true if $\CC$ is the category of
vector spaces over a field $k$, so that $A$ is a $k$-algebra in the
usual sense). In this case, for any free $A$-bimodule $N = A \otimes
V \otimes A$, $V \in \CC$, we have
$$
\Ext^l_{A-\bimod}(N,M)=\Ext^l_{\CC}(V,M)=0, \qquad l \geq 1,
$$
and one can compute the Hochschild cohomology $HH^{\hdot}(A,M)$ by
using the bar-resolution of the $A$-bimodule $\I_A$. This results in
the {\em Hochschild cochain complex} $C^\hdot(A,M)$, where
$$
C^l(A,M) = \Hom_k(A^{\otimes l},M)
$$
for any integer $l \geq 1$. If $M=A$ and $\CC$ is the category of
$k$-vector spaces, one can describe the differential $\delta$ in
this complex as follows: interpret $C^\hdot(A)$ as the graded Lie
algebra of coderivations of the free coalgebra $T^\hdot(A)=
A^{\otimes \hdot}[1]$ generated by the vector space $A$ placed in
degree $-1$ (the bracket $[-,-]$ is given by the graded commutator
of coderivations). Then the multiplication $m:A \otimes A \to A$
gives an element $\delta \in C^2(A)$, and it is easy to check that
$m$ is associative if and only if $\delta^2=[\delta,\delta]=0$. We
assume that this is the case, and the Hochschild differential is
then given by $a \mapsto [\delta,a]$ (for details, and for a
description of the differential for a general $M$, we refer the
reader for instance to \cite[Appendix]{GK}).

\begin{lemma}\label{the.same}
Assume that an $A$-bimodule $M$ is injective as an object in $\CC$,
and that a square-zero extension $\wt{A}$ of a $A$ by $M$ is
identified with $A \oplus M$ as an object in $\CC$. Assume that
under this identification, the multiplication in $\wt{A}$ is
expressed as
\begin{equation}\label{gamma}
a * b = ab + \gamma(a,b)
\end{equation}
for some $\gamma \in C^2(A,M)$. Then $\gamma$ is a Hochschild
cocycle, and it represents the class $\theta_{\wt{A}} \in HH^2(A,M)$.
\end{lemma}

\proof{} The first claim is completely standard; we will prove that
$\theta_{\wt{A}}$ is represented by the cocycle $\gamma$. Fixing an
identification $\wt{A} \cong A \oplus M$ is equivalent to fixing a
map $P_0:A \to \wt{A}$ in $\CC$ which splits the projection $\wt{A}
\to A$. To compute $\theta_{\wt{A}}$, consider the free
$\wt{A}$-$A$-bimodule $\wt{A} \otimes A \otimes A$ which maps
surjectively onto $\I_A \subset A \otimes A$ by $a \otimes b \otimes
c \mapsto ab \otimes c - a \otimes bc$ (this is the first term of
the bar-resolution of the bimodule $\I_A$). The map $\wt{A} \otimes
A \otimes A \to \I_A$ obviously lifts to a map $P:\wt{A} \otimes A
\otimes A \to \Ker T \subset \wt{A} \otimes A$, where $T:\wt{A}
\otimes A \to A$ is the natural projection: we take
\begin{equation*}
P(a \otimes b \otimes c) = (a * P_0(b)) \otimes c - a \otimes bc.
\end{equation*}
To represent $\theta_{\wt{A}}$ by a cocycle, one has to compose $P$
with the projection $\Ker T \to \wt{\I}_A$, notice that it factors
through a map $P':A^{\otimes 3} \to \wt{\I}_A$, then compose $P'$
with the bar-resolution differential $\overline{\delta}:A^{\otimes
4} \to A^{\otimes 3}$, and notice that $P' \circ
\overline{\delta}:A^{\otimes 4} \to \wt{\I}_A$ factors through an
$A$-bimodule map $P'':A^{\otimes 4} \to A$. Explictly,
$\overline{\delta}$ is given by $\overline{\delta}(a \otimes b
\otimes c \otimes d) = ab \otimes c \otimes d - a \otimes bc \otimes
d + a \otimes b \otimes cd$; therefore, since $\gamma$ is a
Hochschild cocycle, we have
\begin{align*}
P''(a \otimes b \otimes c \otimes d)
&= P'(ab \otimes c \otimes d - a \otimes bc \otimes d + a
\otimes b \otimes cd) \\
&= \gamma(ab,c)d - \gamma(a,bc)d + \gamma(a,b)cd\\
&= a\gamma(b,c)d,
\end{align*}
as required.
\endproof

If $M$ is not injective -- for instance, if we want to study
first-order deformations of an algebra $A$ which is not injective as
an object in $\CC$ -- then Lemma~\ref{the.same} no longer
applies. When $\CC$ admits enough injectives, one can circumvent
this problem by replacing $M$ with an injective resolution
$\I^\hdot$. Then as before, the Hochschild cohomology groups
$HH^\hdot(A,M)$ can be computed by the bar-resolution, and this
resulting complex is the complex
$$
C^\hdot(A,\I^\hdot) = \Hom_{\CC}(A^{\otimes\hdot},\I^\hdot)
$$
of Hochschild cochains with values in $\I^\hdot$. On the other hand,
for any square-zero extension $\wt{A}$ of $A$ by $M$, the algebra
$\wt{A}' = (\wt{A} \oplus \I^0)/M$ is a square-zero extension of $A$
by $\I^0$; moreover, we have an exact sequence
\begin{equation}\label{wtdash}
\begin{CD}
0 @>>> \wt{A} @>>> \wt{A}' @>{\tau}>> \I^0/M \subset \I^1.
\end{CD}
\end{equation}
Since $\I^0$ is injective, we in fact can fix an isomorphism
$\wt{A}' \cong A \oplus \I^0$ as objects in $\CC$, so that the
multiplication in $\wt{A}'$ is given by \eqref{gamma} for some
$\gamma \in C^2(A,\I^0)$. Moreover, composing the splitting map
$P_0:A \to \wt{A}'$ with the map $\tau:\wt{A}' \to \I^0/M \subset
\I^1$, we obtain a map $A \to \I^1$, which we can treat as an
element $\gamma' \in C^1(A,\I^1) = \Hom_{\CC}(A,\I^1)$. Then
Lemma~\ref{the.same} can be easily generalized to show that under
these identifications, the Hochschild cohomology class
$\theta_{\wt{A}} \in HH^2(A,M)$ is represented by the cocycle
\begin{equation}\label{gamma-ga}
\gamma + \gamma' \in C^2(A,\I^0) \oplus C^1(A,\I^1)
\end{equation}
in the double complex $C^\hdot(A,\I^\hdot)$. We leave the proof to
the reader.

\section{DG algebras.}

The particular situation where we will use \eqref{gamma-ga} is when
$\CC$ is the tensor category of complexes of vector spaces over the
field $k$, so that algebras in $\CC$ are DG algebras over $k$. We
will only need complexes $K^\hdot$ which are bounded below
($K^p = 0$ for $p \ll 0$). A complex is injective if and only if it
is acyclic. Every complex $K^\hdot$ can be canonically embedded into
the acyclic complex $C(K)^\hdot$, the cone of the identity map
$K^\hdot \to K^\hdot$. The quotient $C(K)^\hdot/K^\hdot$ is by
definition identified with the shifted complex $K^\hdot[1]$. The
construction can be iterated, so that every complex $K^\hdot$ admits
a functorial injective resolution $\I^\hdot$ with $\I^p =
C(K)^\hdot[p]$. Then for any DG algebra $A^\hdot$ and DG-bimodule
$M^\hdot$ over $A^\hdot$, the complex of Hochschild cochains is
given by
\begin{equation}\label{bibar}
C^p(A^\hdot,M^\hdot) = \bigoplus_{0 \leq l \leq p-1}
C^{l,p-l}(A^\hdot,M^\hdot) = \bigoplus_{0 \leq l \leq p-1}
\Hom^l((A^\hdot)^{\otimes p-l},M^\hdot),
\end{equation}
where $\Hom^l$ is the space of vector space maps of degree $l$. The
differential in this complex is the sum of the usual Hochschild
differential which comes from the bar construction, and the
differential which comes from the differentials in the complexes
$A^\hdot$, $M^\hdot$. In the case $M^\hdot=A^\hdot$, one can again
interpret $C^\hdot(A^\hdot)$ as the space of positive-degree
coderivations of the free coassociative coalgebra $T^\hdot(A^\hdot)$
generated by the graded vector space $A^\hdot[1]$; the sum of the
differential $d \in C^1(A^\hdot)$ and the multiplication $m \in
C^2(A^\hdot)$ extends to a coderivation $\delta:T^\hdot(A^\hdot) \to
T^\hdot(A^\hdot)$ of degree $1$ satisfying $\delta^2=0$. Then the
differential in $C^\hdot(A^\hdot)$ is given by $a \mapsto
[\delta,a]$. If we are given a first-order deformation
$\wt{A^\hdot}$ of a DG algebra $A^\hdot$ over $k$, then splitting
the corresponding square-zero extension $\wt{A}'$ in \eqref{wtdash}
is equivalent to fixing an isomorphism $\wt{A^\hdot} \cong
A^\hdot\langle \eps \rangle$ of graded vector spaces. Then the
multiplication and the differential $\wt{d}$ in $\wt{A^\hdot}$ are
given by
\begin{align*}
a * b &= ab + \gamma_2(a,b)\eps,\\
\wt{d}(a) &= d(a) + \gamma_1(a)\eps
\end{align*}
for some $\gamma_1 \in C^{1,1}(A^\hdot)$, $\gamma_2 \in C^{2,0}(A)$,
and by \eqref{gamma-ga}, $\gamma = \gamma_1 + \gamma_2$ is a
Hochschild cocycle representing the class $\theta_{\wt{A}} \in
HH^2(A)$.

However, we will also need a variation of the Hochschild cohomology
construction specific to DG algebras. Namely, given a flat DG
algebra $A^\hdot$ in some tensor category $\CC$, one can invert
quasiisomorphism in the category of $A^\hdot$-bimodules and obtain
the derived category $\D(A^\hdot)$ of DG $A^\hdot$-bimodules. This
gives a DG version
$$
HH^\hdot_\D(A^\hdot) = \RHom^\hdot_{\D(A^\hdot)}(\I_{A^\hdot},A)
$$
of the Hochschild cohomology. We have a canonical map
\begin{equation}\label{to.d}
HH^\hdot(A^\hdot) \to HH^\hdot_\D(A^\hdot).
\end{equation}
If $\C$ is the category of $k$-vector spaces, the groups
$HH^\hdot_\D(A^\hdot)$ can be computed by the same complex
\eqref{bibar}, but without the condition $l \geq 0$.

The point of introducing the groups $HH^\hdot_\D(A^\hdot)$ is that
they control deformations ``up to quasiisomorphism''. Namely, recall
that the category of DG algebras up to a quasiisomorphism is
obtained from the category of DG algebras by formally inverting all
algebra maps which are quasiisomorphisms -- in other words, DG
algebras $A^\hdot$, $B^\hdot$ are quasiisomorphic if there exists a
chain of quasiisomorphisms 
\begin{equation}\label{qs}
\begin{CD}
A^\hdot @<<< A_1^\hdot @>>> A_2^\hdot @<<< \dots @<<< A^\hdot_n @>>>
B^\hdot.
\end{CD}
\end{equation}
Then we have the following fact.

\begin{lemma}\label{bootstrap}
Assume that $\CC$ is the category of vector spaces of a field $k$ of
characteristic $\cchar k = 0$. Assume given a DG algebra $A^\hdot$
in $\CC$ and its formal deformation $\wt{A^\hdot}$. Then for any
integer $p \geq 1$, the DG algebra $\wt{A^\hdot}/h^{p+1}$ is
quasiisomorphic to the DG algebra $A^\hdot[h]/h^{p+1}$ if and only
if the Kodaira-Spencer class $\Theta_{\wt{A^\hdot}}$ of the
deformation $\wt{A^\hdot}$ vanishes after projection to
$HH^2_\D(\wt{A^\hdot})/h^p$.
\end{lemma}

\proof{} This is, in a sense, a DG version of Z. Ran's $T_1$-lifting
principle \cite{ran}. To control quasiisomorphisms, it is convenient
to use the notion of an $A_\infty$-morphism. Recall (see
e.g. \cite{kell}) that an $A_\infty$-morphism $\iota$ between DG
algebras $A$ and $B$ is by definition a map $\iota:T^\hdot(A^\hdot)
\to T^\hdot(B^\hdot)$ between the free coassociative DG coalgebras
$T^\hdot(A^\hdot)$, $T^\hdot(B^\hdot)$ generated by $A[1]$ and
$B[1]$ such that $\delta_B \circ \iota = \iota \circ \delta_A$,
where $\delta_A$ and $\delta_B$ are Hochschild differentials on
$T^\hdot(A^\hdot)$, $T^\hdot(B^\hdot)$. Every DG algebra map
obviously induces an $A_\infty$-morphism. However, if a DG algebra
map is a quasiisomorphism, then the corresponding $A_\infty$-map
$\iota$ is {\em invertible} (that is, there exists an $A_\infty$-map
$\iota^{-1}$ such that $\iota \circ \iota^{-1}$ and $\iota^{-1}
\circ \iota$ are identical on cohomology). Therefore if two DG
algebras $A^\hdot$, $B^\hdot$ are quasiisomorphic, not only there
exists a chain of quasiisomorphisms \eqref{qs} -- there in fact
exists a single $A_\infty$-quasiisomorphism $\iota:T^\hdot(A^\hdot)
\to T^\hdot(B^\hdot)$.

Now, fix a graded vector space isomorphism $\wt{A^\hdot} \cong
A^\hdot \otimes_k k[h]$; under this isomorphism, the multiplication
and the differential in $\wt{A}$ are given by
\begin{equation}\label{a_infty}
\begin{aligned}
d &= d_0 + \sum_{l \geq 1}h^ld_l \\
m &= m_0 + \sum_{l \geq 1}h^lm_l,
\end{aligned}
\end{equation}
where $m_0$ and $d_0$ are the multiplication and the differential in
$A^\hdot$. The Hochschild cohomology $HH^\hdot_\D(\wt{A^\hdot})$ of
the DG $k[h]$-algebra $\wt{A^\hdot}$ can be computed by the same
Hochschild complex $C^\hdot(A^\hdot)[h]$ as $HH^\hdot_\D(A^\hdot
\otimes_k k[h])$, but with a different differential: the
differentials $\delta$ and $\delta_0$ computing
$HH^\hdot_\D(\wt{A^\hdot})$ and $HH^\hdot_\D(A^\hdot)[h]$ are given
by $\delta(a) = [d+m,a]$, $\delta_0(a) = [d_0+m_0,a]$. The quotient
$\wt{A^\hdot}/h^{p+1}$ is quasiisomorphic to $A^\hdot[h]/h^{p+1}$ if
and only if there exists a $k[h]/h^{p+1}$-linear
$A_\infty$-quasiisomorphism between them -- in other words, a
coalgebra map $\iota_{p+1}: T^\hdot(A^\hdot)[h] \to
T^\hdot(A^\hdot)[h]$ such that $\delta \circ \iota_{p+1} =
\iota_{p+1} \circ \delta_0 \mod h^{p+1}$.

To compute the cocycle $\Theta=\Theta_{\wt{A^\hdot}}$, we can use
\eqref{gamma-ga}. Namely, we replace $h$ with $h+\eps$ in
\eqref{a_infty}, and we conclude that the image of the
Kodaira-Spencer class $\Theta$ in the group $H^2_\D(\wt{A^\hdot})$
is represented by the cocycle
$$
Q=\sum_{l \geq 1}lh^l(m_l+d_l) \in C^\hdot(A^\hdot)[h]
$$
of total degree $2$. To prove the claim of the Lemma, use induction
on $p$. Assume by induction that $\Theta = 0 \mod h^p$ and that
there exists a map $\iota_p:T^\hdot(A^\hdot)[h] \to
T^\hdot(A^\hdot)[h]$ such that $\delta \circ \iota_p = \iota_p \circ
\delta_0 \mod h^p$.  Then $\iota_p(Q)$ is divisible by $h^p$, and it
represents the class $\Theta$. Thus $\Theta=0 \mod h^{p+1}$ if and
only if $\iota_p(Q)=h^p\delta(\gamma) \mod h^{p+1}$ for some
Hochschild cochain $\gamma \in C^\hdot(A^\hdot)$ of total degree
$1$. Since $\delta=\delta_0 \mod h$, this can be rewritten as
$$
\iota_p(Q)=h^p[\delta,\gamma]=h^p[\delta_0,\gamma] \mod h^{p+1},
$$
which is in turn equivalent to $\delta \circ \iota_{p+1} =
\iota_{p+1} \circ \delta_0 \mod h^{p+1}$, where we set $\iota_{p+1} =
\iota_p + \frac{1}{p}h^p\gamma$. This proves the claim.
\endproof

\begin{remark}
The point in the above proof where we do need to consider
$A_\infty$-morphisms is in the construction of the correction term
$\gamma$: the cochain $\gamma \in C^2(A^\hdot)$ in the Hochschild
complex \eqref{bibar} can have non-trivial components in
$C^{l,2-l}(A^\hdot)$ with $l < 0$.
\end{remark}

When $\CC$ is a general tensor category, the relation between the
Kodaira-Spencer class and the triviality of deformations of DG
algebras is a difficult subject better left untouched in the present
paper. However, the Kodaira-Spencer class itself is perfectly well
defined.

\section{Obstructions to formality.}
We can now proceed to our objective -- the study of
formality. Assume given a DG algebra $A^\hdot$ in a tensor abelian
category $\CC$. We want to study whether $A^\hdot$ is formal -- that
is, quasiisomorphic to the cohomology DG algebra $H^\hdot(A^\hdot)$
(with trivial differential). To do this, consider the canonical
filtration $F_\idot$ on $A^\hdot$ -- that is, set
$$
F_k A^p =
\begin{cases}
0, &\qquad \text{ if }p>k,\\
\Ker d:A^p \to A^{p+1}, &\qquad \text{ if }p=k,\\
A^p, &\qquad \text{ if }p<k,
\end{cases}
$$ 
and denote by $B^\hdot \cong \gr^F A^\hdot$ the associated graded
quotient. The canonical filtration induces a filtration $F_\idot$ on
the Hochschild cohomology complex $HH^\hdot_\D(A^\hdot)$; for any
two integers $p \leq q$ we denote
$$
HH^\hdot_\D(A^\hdot)_{(p,q)} =
F_pHH^\hdot_\D(A^\hdot)/F_qHH^\hdot_\D(A^\hdot).
$$
The associated graded quotient $\gr^F_\idot HH^\hdot_\D(A^\hdot)$ is
naturally quasiisomorphic to $HH^\hdot_\D(B^\hdot)$. Denote the
induced grading on the complex $HH^\hdot_\D(B^\hdot)$ by
$HH^\hdot_\D(B^\hdot)_\idot$, so that we have
\begin{equation}\label{gr}
\gr^p_F HH^\hdot_\D(A^\hdot) = HH^\hdot_\D(A^\hdot)_{(p,p)} \cong
HH^\hdot_\D(B^\hdot)_p
\end{equation}
for any integer $p$. If $\CC$ is the category of vector spaces over
$k$, one can compute $HH^\hdot_{\D}(B^\hdot)$ by means of the
bar-construction; then the complex computing the component
$HH^\hdot_{\D}(B^\hdot)_p$ of degree $p$ consists of spaces
\begin{equation}\label{bar}
\Hom^{q-p}((B^\hdot)^{\otimes q},B^\hdot).
\end{equation}

Now, the natural quasiisomorphism $\gr^F A^\hdot \to H^\hdot(A)$ is
compatible with the multiplication. Therefore the question of
formality of the DG algebra $A^\hdot$ is equivalent to the existence
of a quasiisomorphism between the DG algebras $A^\hdot$ and $B^\hdot
= \gr^F A^\hdot$. To refine this, consider the Rees algebra
$$
\wt{A^\hdot} = \bigoplus_p F_pA^\hdot;
$$
this is a graded DG-algebra in $\CC[h]$, with $h$ of degree $1$
given by the natural embedding $h:F_\idot A^\hdot \to
F_{\idot+1}A^\hdot$. Then $\wt{A^\hdot}/h\wt{A^\hdot} \cong
B^\hdot$, so that $\wt{A^\hdot}$ is a formal deformation of the DG
algebra $B^\hdot$.

\begin{defn}
The DG algebra $A^\hdot$ is {\em $p$-formal} for some integer $p
\geq 1$ if $\wt{A^\hdot}/h^{p+1}$ is quasiisomorphic to
$B^\hdot[h]/h^{p+1}$.
\end{defn}

\begin{lemma}
A DG algebra $A^\hdot$ in the category of vector spaces over a field
$k$ is formal if and only if it is $p$-formal for every $p \geq 1$.
\end{lemma}

\proof{} The ``only if'' part is obvious: if $A^\hdot$ is formal,
then its deformation $\wt{A^\hdot}$ is trivial, so that for every $p
\geq 1$, the truncation $\wt{A^\hdot}/h^{p+1}$ is quasiisomorphic to
$B^\hdot[h]/h^{p+1}$. Conversely, if $A^\hdot$ is $p$-formal, then
we can choose an $A_\infty$-map $s:B^\hdot \to \wt{A^\hdot}/h^{p+1}$
which spltis the natural projection $\wt{A^\hdot}/h^{p+1} \to
B^\hdot$. Moreover, if $A^\hdot$ is $p$-formal for every $p$, then
we can choose these splitting maps in a compatible way and obtain an
$A_\infty$-map $s:B^\hdot \to \wh{A^\hdot}$ from $B^\hdot$ to the
completion $\wh{A^\hdot} = \lim_{\gets}\wt{A^\hdot}/h^{p+1}$ of the
algebra $\wt{A^\hdot}$ with respect to the $h$-adic
topology. Explicitly, this completion $\wh{A^\hdot}$ is a graded DG
algebra, whose component of degree $m$ is equal to
$$
\wh{A^\hdot}_m = \lim_{l \to -\infty}F_mA^\hdot/F_lA^\hdot.
$$
However, for any complex $V^\hdot$ with the canonical filtration
$F_\idot V^\hdot$, the inverse system $V^\hdot/F_lV^\hdot$ obviously
stabilizes in any degree at a finite step ($V^m/F_lV^m$ stops
depending on $l$ when $l<m-1$). Therefore the completion
$\wh{A^\hdot}$ is actually isomorphic to $\wt{A^\hdot}$, so that we
have an $A_\infty$-map $B^\hdot \to \wt{A^\hdot}$ which splits the
projection $\wt{A^\hdot} \to B^\hdot$. Evaluating at $h=1$, we get
an $A_\infty$-quasiisomorphism $B^\hdot \cong A^\hdot$.
\endproof

To measure $p$-formality for all $p \geq $, we use
Lemma~\ref{bootstrap}. Denote by
$$
Q_{A^\hdot} \in HH^2_\D(\wt{A^\hdot})
$$
the image of the Kodaira-Spencer class $\Theta_{\wt{A^\hdot}}$ under
the natural map \eqref{to.d}. The grading on the Rees algebra
induces a natural grading on $HH^\hdot_\D(\wt{A^\hdot})$, with the
component of degree $p$ canonically identified with
$F_pHH^\hdot_\D(A^\hdot)$. Since $h$ is of degree $1$, the class
$Q_{A^\hdot}$ is of degree $-1$, so that in fact we have the
canonical class
$$
Q_{A^\hdot} \in F_{-1}HH^2_{\D}(A^\hdot).
$$
This class is a version of the so-called {\em higher Massey
products} in the DG algebra $A^\hdot$ -- all of them in one
package. Modulo $h^p$, this class restricts to a class in the degree
$-1$-component of $HH^2_\D(\wt{A^\hdot})/h^p$, which is identified
with $HH^2_\D(A^\hdot)_{(-p,-1)}$.

\medskip

Let now $\CC$ be the category of sheaves of $\calo_X$-modules on a
scheme $X$ over a field of characteristic $0$. Then in addition to
tensor structure, the category $\CC = \Coh(X)$ has inner $\Hom$ and
its derived functors, which we denote by $\Rhom^\hdot$. The same is
true for $C^\hdot(\CC)$ and for the category of $\A^\hdot$-bimodules
for some DG algebra $\A^\hdot$ in $\CC$. This allows to refine the
construction of Hochschild cohomology: we can define the Hochschild
cohomology complex
$$
\hh^\hdot_{\D}(\A^\hdot) =
\Rhom^\hdot_{\D(\A^\hdot)}(\I_{\A^\hdot},\A^\hdot)
$$
of sheaves of $\calo_X$-modules on $X$, and we have
$HH^\hdot_{\D}(\A^\hdot) \cong
\HH^\hdot(X,\hh^\hdot_{\D}(\A^\hdot))$. The Rees algebra, the
canonical filtration on $\A^\hdot$, and $\hh^\hdot_{\D}(\A^\hdot)$
are also well-defined on the level of sheaves, and so is the grading
\eqref{gr}.

We note that if $X$ is Noetherian, then the inner $\Hom$ between two
{\em coherent} sheaves of $\calo_X$-modules is also coherent, and
$\Rhom^\hdot$ is a complex with coherent homology sheaves.

\begin{theorem}\label{Q}
Let $\A^\hdot$ be a DG algebra of flat sheaves of $\calo_X$-modules
on a Noetherian reduced irreducible scheme $X$ over a field of
characteristic $0$. Let $\B^\hdot$ be the homology algebra of the DG
algebra $\A^\hdot$. Assume that the sheaves $\B^\hdot$ are flat and
coherent on $X$. Assume also that for every integers $l$, $i$, the
component $\hh^i_{\D}(\B^\hdot)_{-l}$ of degree $(-l)$ of $i$-th
Hochschild cohomology sheaf $\hh^i_{\D}(\B^\hdot)$ with respect to
the grading \eqref{gr} is coherent and flat.
\begin{enumerate}
\item If the fiber $\A^\hdot_x$ is formal for a generic point $x \in
X$, then it is formal for an arbitrary point $x \in X$.
\item Assume in addition that for every integer $l \geq 1$, we have
$$
H^0(X,(\hh^2_{\D}(\B^\hdot))_{-l})=0.
$$
Then the DG algebra $\A^\hdot_x$ is formal for every point $x \in
X$.
\end{enumerate}
\end{theorem}

\proof{} Use induction and Lemma~\ref{bootstrap}. Assume that the
fiber $\A^\hdot_x$ is $p$-formal for some integer $p \geq 1$ and
every point $x \in X$. Consider the spectral sequence in the
category of sheaves of $\calo_X$-modules on $X$ associated to the
filtration $F_\idot$ on the complex
$\hh^\hdot_\D(\A^\hdot)_{-p,-1}$. In terms of the Rees algebra, this
is the $h$-adic filtration on the $(-1)$-component of
$\hh^\hdot_\D(\wt{A^\hdot}/h^p)$. The terms of this spectral
sequence are the Hochschild cohomology sheaves
$\hh^\hdot_{\D}(\B^\hdot)_{-l}$, $1 \leq l \leq p$, and by
assumption, these are flat coherent sheaves on $X$. The fiber of the
spectral sequence at any point $x \in X$ gives the corresponding
spectral sequence for the fiber $\A^\hdot_x$. Since by assumption
the algebra $\A^\hdot_x$ is $p$-formal for every $x \in X$, the
differential in the spectral sequence vanishes at every
point. Therefore the spectral sequence degenerates, and the homology
sheaves of the complex $\hh^\hdot_\D(\A^\hdot)_{-p,-1}$ are iterated
extensions of the homology sheaves of the complexes
$\hh^\hdot_{\D}(\B^\hdot)_{-l}$, $1 \leq l \leq p$. In particular,
the second homology sheaf $\hh^2_\D(\A^\hdot)_{-p,-1}$ is flat and
coherent. Thus if the reduction $\overline{Q}_{\A^\hdot} =
Q_{\A^\hdot} \mod h^p \in H^0(X,\hh^2_\D(\A^\hdot)_{-p,-1})$
vanishes at the generic point of $X$, it vanishes everywhere, and
$\A^\hdot$ is $(p+1)$-formal everywhere by
Lemma~\ref{bootstrap}. This proves \thetag{i}. Moreover, in the
assumptions of \thetag{ii} we in fact have
$H^0(X,(\hh^2_{\D}(\A^\hdot))_{-p,-1})=0$, which by
Lemma~\ref{bootstrap} again proves that $\A^\hdot$ is $(p+1)$-formal
at every point $x \in X$.
\endproof

\bigskip

\noindent
{\sc Steklov Math Institute\\
Moscow, USSR}

\bigskip

\noindent
{\em E-mail address\/}: {\tt kaledin@mccme.ru}

\end{document}